\documentclass[12pt]{article}

\usepackage{amssymb,amsmath,amsthm,amsfonts,amscd,bbm}
\usepackage{enumitem}
\usepackage{pdflscape}
\usepackage[paper=portrait,pagesize]{typearea}
\usepackage{caption}
\usepackage{bm}

\usepackage{geometry} 

\usepackage{ifpdf}
\ifpdf
\usepackage[pdftex]{graphicx}
\else
\usepackage[dvips]{graphicx}
\fi
\usepackage[all]{xy}
\usepackage{hyperref}
\usepackage{xcolor}
\usepackage{graphicx}
\usepackage[outdir=./]{epstopdf}
\usepackage{enumitem}
\usepackage{tensor}
\usepackage{tikz-cd}
\usepackage{multicol}
\usepackage{mathtools}
\usepackage{tocvsec2}
\usepackage{bbm}
\usepackage{longtable}
\usepackage{fullpage}

\newcommand{\googlebooks}[1]{(preview at \href{https://books.google.com/books?id=#1}{google books})}

\newcommand{\numdam}[1]{}
\usepackage{mathrsfs}

\usepackage{mathrsfs}
\DeclareMathAlphabet{\mathpzc}{OT1}{pzc}{m}{it}

\def\semicolon{;}
\def\applytolist#1{
    \expandafter\def\csname multi#1\endcsname##1{
        \def\multiack{##1}\ifx\multiack\semicolon
            \def\next{\relax}
        \else
            \csname #1\endcsname{##1}
            \def\next{\csname multi#1\endcsname}
        \fi
        \next}
    \csname multi#1\endcsname}

\def\calc#1{\expandafter\def\csname c#1\endcsname{{\mathcal #1}}}
\applytolist{calc}QWERTYUIOPLKJHGFDSAZXCVBNM;
\def\bbc#1{\expandafter\def\csname bb#1\endcsname{{\mathbb #1}}}
\applytolist{bbc}QWERTYUIOPLKJHGFDSAZXCVBNM;
\def\bfc#1{\expandafter\def\csname bf#1\endcsname{{\mathbf #1}}}
\applytolist{bfc}QWERTYUIOPLKJHGFDSAZXCVBNM;
\def\sfc#1{\expandafter\def\csname s#1\endcsname{{\sf #1}}}
\applytolist{sfc}QWERTYUIOPLKJHGFDSAZXCVBNM;
\def\fc#1{\expandafter\def\csname f#1\endcsname{{\mathfrak #1}}}
\applytolist{fc}QWERTYUIOPLKJHGFDSAZXCVBNM;

\usepackage{tikz}
\usepackage{tikz-cd}

\makeatletter
\def\fixtikzforbreqn#1#2{%
  \protected\edef#1{\noexpand\ifmmode\mathchar\the\mathcode`#2 \noexpand\else#2\noexpand\fi}%
}
\fixtikzforbreqn\tikz@nonactivesemicolon;
\fixtikzforbreqn\tikz@nonactivecolon:
\fixtikzforbreqn\tikz@nonactivebar|
\fixtikzforbreqn\tikz@nonactiveexlmark!
\makeatother

\usepackage{breqn}   

\usetikzlibrary{arrows,backgrounds,patterns.meta}
\usetikzlibrary{positioning,shadings,cd}
\usetikzlibrary{shapes}
\usetikzlibrary{backgrounds}
\usetikzlibrary{decorations,decorations.pathreplacing,decorations.markings,decorations.pathmorphing}
\usetikzlibrary{fit,calc,through}
\usetikzlibrary{external}
\usetikzlibrary{arrows}
\tikzset{vertex/.style = {shape=circle,draw,fill=black,inner sep=0pt,minimum size=5pt}}
\tikzset{edge/.style = {->,> = latex', bend right}}
\tikzset{
	super thick/.style={line width=3pt}
}
\tikzset{
    quadruple/.style args={[#1] in [#2] in [#3] in [#4]}{
        #1,preaction={preaction={preaction={draw,#4},draw,#3}, draw,#2}
    }
}
\tikzstyle{shaded}=[fill=red!10!blue!20!gray!30!white]
\tikzstyle{unshaded}=[fill=white]
\tikzstyle{empty box}=[circle, draw, thick, fill=white, opaque, inner sep=2mm]
\tikzstyle{annular}=[scale=.7, inner sep=1mm, baseline]
\tikzstyle{rectangular}=[scale=.75, inner sep=1mm, baseline=-.1cm]
\tikzstyle{mid>}=[decoration={markings, mark=at position 0.5 with {\arrow{>}}}, postaction={decorate}]
\tikzstyle{mid<}=[decoration={markings, mark=at position 0.5 with {\arrow{<}}}, postaction={decorate}]
\tikzstyle{over}=[double, draw=white, super thick, double=]
\tikzstyle{snake}=[decorate, decoration={snake, segment length=1mm, amplitude=.3mm}]
\tikzstyle{saw}=[decorate, decoration={saw, segment length=.7mm, amplitude=.25mm}]

\tikzstyle{coupon}=[draw, very thick, rectangle, rounded corners=5pt]
\tikzset{Rightarrow/.style={double equal sign distance,>={Implies},->},
triplecd/.style={-,preaction={draw,Rightarrow}},
quadruplecd/.style={preaction={draw,Rightarrow,
shorten >=0pt
},
shorten >=1pt,
-,double,double
distance=0.2pt}}
\tikzset{
    tripleline/.style args={[#1] in [#2] in [#3]}{
        #1,preaction={preaction={draw,#3},draw,#2}
    }
}
\tikzstyle{triple}=[tripleline={[line width=.15mm,black] in
      [line width=.7mm,white] in
      [line width=1mm,black]}] 
\tikzset{
    quadrupleline/.style args={[#1] in [#2] in [#3] in [#4]}{
        #1,preaction={preaction={preaction={draw,#4},draw,#3}, draw,#2}
    }
}
\tikzstyle{quadruple}=[quadrupleline={[line width=.3mm,white] in
      [line width=.6mm,black] in
      [line width=1.2mm,white] in
      [line width=1.5mm,black]}]

\makeatletter
\newcommand{\xMapsto}[2][]{\ext@arrow 0599{\Mapstofill@}{#1}{#2}}
\def\Mapstofill@{\arrowfill@{\Mapstochar\Relbar}\Relbar\Rightarrow}
\makeatother

\theoremstyle{plain}
\newtheorem{thm}{Theorem}[section]
\newtheorem*{thm*}{Theorem}

\newtheorem*{cor*}{Corollary}

\newtheorem*{conj*}{Conjecture}
\newtheorem{lem}[thm]{Lemma}
\newtheorem*{lem*}{Lemma}

\newtheorem*{quest*}{Question}
\newtheorem*{claim*}{Claim}

\theoremstyle{definition}
\newtheorem{defn}[thm]{Definition}

\newtheorem{sub-ex}[thm]{Sub-Example}
\newtheorem{counter-ex}[thm]{Counter-Example}
\newtheorem{rem}[thm]{Remark}
\newtheorem*{rem*}{Remark}

\usepackage{xcolor}
\definecolor{dark-red}{rgb}{0.7,0.25,0.25}
\definecolor{dark-blue}{rgb}{0.15,0.15,0.55}
\definecolor{medium-blue}{rgb}{0,0,.8}
\definecolor{DarkGreen}{RGB}{0,150,0}
\definecolor{rho}{named}{red}
\hypersetup{
   colorlinks, linkcolor={purple},
   citecolor={medium-blue}, urlcolor={medium-blue}
}

\newcommand{\id}{\operatorname{id}}

\newcommand{\rad}{\operatorname{rad}}

\newcommand{\Mod}{{\sf Mod}}

\newcommand{\Bim}{{\sf Bim}}

\newcommand{\QuivCon}{{\sf QuivCon}}
\newcommand{\PathAlg}{{\sf PathAlg}}

\newcommand{\BdQuivCon}{{\sf BdQuivCon}}
\newcommand{\BSA}{{\bf BSA}}
\def\altdb{\vadjust{\vbox to 0pt{\vss\hbox{\kern \hsize
\quad{\dbend}}\kern\baselineskip\kern-10pt}}}

\newcommand{\noshow}[1]{}
\title{Quiver connections and bimodules of basic algebras}
\author{Sean Thompson}
\date{}

\begin{document}

\maketitle

\begin{abstract}
Motivated by the problem of classifying quantum symmetries of non-semisimple, finite-dimensional associative algebras, we define a notion of connection between bounded quivers and build a bicategory of bounded quivers and quiver connections. We prove this bicategory is equivalent to a bicategory of basic algebras,  bimodules, and intertwiners with some additional structure. 
\end{abstract}

\tableofcontents

\section{Introduction}
Basic algebras are an important class of finite dimensional algebras. An algebra $A$ with a complete set of primitive orthogonal idempotents $\{e_1,\hdots\,e_n\}$ is basic if $e_iA\not\cong e_jA$ for all $i\neq j$. Since every finite dimensional algebra is Morita equivalent to a basic algebra \cite{morita}, these algebras play a crucial role in understanding the representation theory of arbitrary finite dimensional algebras. 

Quivers are useful combinatorial tool for understanding basic algebras and their representations. A well-known theorem of Gabriel \cite{quiver} characterizes basic algebras in terms of quivers. This theorem states that any finite dimensional, basic algebra $A$ with a complete set of $n$ primitive orthogonal idempotents is the quotient of the path algebra of a quiver $Q_A$ on $n$ vertices by an admissible ideal $I_A$. The data $(Q_A,I_A)$, commonly called a \textit{bound quiver}, is not uniquely determined by $A$ alone. The possible choices of bound quivers realizing $A$ are parameterized by what we will call \textit{quiver data}.

Quiver data for a basic algebra $A$ consists of linear lifting maps $\delta^1_A:A/\rad A\to A$ and $\delta^2_A:\rad A/\rad^2 A\to \rad A$ that satisfy some coherence conditions (Definition \ref{quivdat}).  It is useful to think of $\delta^1_A$ as assigning each vertex in $Q_A$ to a unique primitive idempotent $e_i\in A$. Similarly, the choice of $\delta^2_A$ is analogous to assigning edges in $Q_A$ to basis elements of the space $\rad A/\rad^2 A$. Since a quiver is defined by its sets of vertices and edges, these data are sufficient to fully parameterize the map from $A$ to $Q_A/I_A$.

Recently, there has been significant interest in understanding \textit{actions} of fusion categories on algebras, characterized by a linear monoidal functor from an abstract fusion category to $\text{Bim}(A)$, where $A$ is some associative algebra. In finite dimensions, this has been completely characterized in cases when $A$ is semisimple, but the case for $A$ non-semisimple is poorly understood (see section \ref{application}). To help us better understand the category of bimodules of a basic algebra, a natural question to ask is: can we extend extend the above story relating basic algebras and quivers to bimodules?

Motivated by Ocneanu's theory of biunitary connections and recent generalizations \cite{connect2, connect1, fusion}, we define a bicategory of bound quivers and quiver ``connections'', called \BdQuivCon, where a quiver connection can be understood to be edges pointing from one quiver into another, along with a map to move paths across the connection (Definition \ref{connect}).

We introduce a definition of a bimodule version of a basic algebra which, given two basic algebras $A$ and $B$, consists of a dualizable $A-B$ bimodule $M$ such that $\rad A M\cong M\rad B$, together with \textit{bimodule quiver data}. These will be linear liftings $\delta^1_M:M/\rad M\to M$ and $\delta^2_M:\rad M/\rad^2 M\to \rad M$ that satisfy similar coherence conditions to those for the algebra quiver data (Definition \ref{bimquivdat}). These assemble into a bicategory $\BSA$ of basic algebras (Definition \ref{bsa}). We then prove the following theorem: 

\begin{thm*}[1]\label{ThmA}
$\BdQuivCon$ and $\BSA$ are equivalent as 2-categories.    
\end{thm*}
The above theorem can be viewed as a 2-categorical version of Gabriel's original result on bound-quivers.

\subsection{Application: actions of fusion categories on truncated path algebras}\label{application}
One of our original motivations is to study actions of fusion categories on non-semisimple algebras. Recall a \textit{fusion category} (over $\mathbbm{C}$) is a finitely semisimple rigid tensor category, with simple unit object \cite{fusion2}. An \textit{action} of a fusion category $\mathcal{C}$ on an associative algebra $A$ is a linear monoidal functor $\mathcal{C}\rightarrow \text{Bim}(A)$.

There is a sense in which we understand all possible actions of fusion categories on a semisimple algebra, since a semisimple algebra is the same data as a semisimple category with a choice of a generating object. An action of $\cC$ on $A$ is then is the same as a module category structure on $\Mod(A)$, using semisimplicity and the Eilenberg-Watts theorem \cite{eilen}\cite{watts}. In particular, for $F:\cC\to\Bim(A)$ and $X\in\cC$, $F(X)\otimes_A\cdot$ equips $\Mod(A)$ with the structure of a module category, and conversely. Semisimple module categories of a fusion category can in turn be understood in terms of Morita classes of algebras internal to $\mathcal{C}$ \cite{module}.
However, to our knowledge the case of finite-dimensional, non-semisimple algebras is largely unexplored territory. One motivation for this paper arises from the goal of taking actions of fusion categories on $\mathbbm{k}^n$ and studying ``liftings'' of the action on some ``de-semisimplifcation" $A$. Our construction of the category $\BSA$ is designed to make precise the notion of lifting a $\mathbbm{k}^n$ bimodule to an $A$ bimodule, which we plan to expand further in future work. However, we have the following definition:

\begin{defn}
    A \textit{basic} action of a fusion category $\mathcal{C}$ on a basic algebra $A$ with quiver data $(\delta^1_A, \delta^2_A)$ is a linear monoidal functor $\mathcal{C}\rightarrow \BSA(A,\delta^1_A, \delta^2_A)$.
\end{defn}

As a corollary of Theorem \ref{ThmA} basic actions are parameterized by linear monoidal functors $\mathcal{C}\rightarrow \BdQuivCon(Q,I)$, where $(Q,I)$ are determined by $A$. This translates the problem into a more combinatorial setting. 

As an example of the utility of this approach, we will consider basic actions of fusion categories on \textit{truncated path algebras}. A truncated path algebra for a quiver $Q$ is simply $Q/I_{n}$, where $I_{n}$ is the ideal spanned by all paths of length at least $n$. For $n\ge 2$ this is admissible, so the truncated path algebras are basic algebras. We consider these as objects in $\BSA$ with the obvious quiver data. If $Q$ is a quiver, denote by $V(Q)$ the set of vertices and $\text{Vec}(V(Q))$ the semisimple linear category of vector bundles on $V(Q)$. If $\mathcal{M}$ is a $\mathcal{C}$-module category, let $\mathcal{C}^{*}_{\mathcal{M}}$ denote the \textit{dual category} in the sense of \cite{tensor}. We have the following theorem:

\begin{thm*}[2]\label{ThmB} Let $\mathcal{C}$ be a fusion category, and $Q/I_{n}$ be a truncated path algebra. Then isomorphism classes of basic actions of $\mathcal{C}$ on $Q/I_{n}$ are parameterized by isomorphism classes of module category structures on the semisimple category $\text{Vec}(V(Q))$ together with a class of object in the dual multi-fusion category  $F\in \mathcal{C}^{*}_{\text{Vec}(V(Q))}$ whose fusion graph is isomorphic to $Q$.
\end{thm*}

\begin{proof}

First note that for any two quivers $P$ and $Q$, the natural inclusion

$$\BdQuivCon((P,I_{n}),(Q,I_{n}))\hookrightarrow \QuivCon(P,Q)$$

is fully faithful, and thus basic actions on \textit{any} truncated path algebra are simply paramtered by linear monoidal functors $F:\mathcal{C}\rightarrow \QuivCon(Q)$. But this is the same data as a module category structure on $\text{Vec}(V(Q))$ together with the structure of a $\mathcal{C}$-module functor on the the endofunctor on $\text{Vec}(V(Q))$ corresponding to $Q$.
\end{proof}

We encourage the reader to compare the above result with \cite{fusion,connect1} and also \cite{finite}, which has tensor algebras in the fusion category $\mathcal{C}$ parameterized by the same data.

\section{Connection Categories}

Recall that a \textit{quiver} $Q$ is an oriented multi-graph with finite vertex and edge sets. Here we introduce the following definition of a connection between two quivers. This will serve as a 1-morphism in a 2-category of quivers.

\begin{defn}\label{connect}
    For quivers $G$ and $H$, a \textbf{connection} is constructed as follows. Fix a field $\mathbbm{k}$, and define $E^{G}_{g,g'}=\mathbbm{k}[E(G)(g\to g')], E^{H}_{h,h'}=\mathbbm{k}[E(H)(h\to h')]$. Then let $\Gamma$ be a family of finite dimensional vector spaces $\{\Gamma_{g,h}\}_{g\in V(G),h\in V(H)}$, and $U$ be a family of linear isomorphisms:
    \[U_{g,h}:\bigoplus_{g'\in V(G)} E^{G}_{g,g'}\otimes \Gamma_{g',h}\to \bigoplus_{h'\in V(H)} \Gamma_{g,h'}\otimes E^H_{h',h}.\] 
    The pair $(U,\Gamma)$ is a $G,H$ connection, and it is useful think of $\Gamma$ as a quiver with edges from vertices in $G$ pointing into vertices in $H$. There is an isomorphism between length one edges and basis vectors of vector spaces, so throughout the paper we use paths and vector spaces interchangeably.
\end{defn}
\begin{defn}
    Now we define a 2-category $\QuivCon$ whose

\begin{itemize}
    \item Objects are quivers
    \item 1-morphisms are connections, with horizontal composition defined as follows: if $(U,\Gamma)$ is a $G,H$ connection and $(V,\Delta)$ is a $H,K$ connection, then
    \[(V,\Delta)\otimes (U,\Gamma):= (\Gamma\otimes\Delta), \left(\id_{\Gamma}\otimes V \right)  \circ \left(U\otimes \id_{\Delta},\right)\]
    where \[(\Gamma\otimes \Delta)_{i,j}=\bigoplus_{k\in V(H)} \Gamma_{i,k}\otimes \Delta_{k,j}\]
    \item A 2-morphism between two $G,H$ connections $(U,\Gamma)$ and $(V,\Delta)$ $f$ is a family of linear maps $\{f_{g,h}:\Gamma_{g,h}\to\Delta_{g,h}\}_{g\in V(G), h\in V(H)}$ satisfying
    \[V_{g,g',h,h'}\circ \left(\id_{E^G_{g,g'}}\otimes f_{g',h}\right)=\left(f_{g,h'}\otimes \id_{E^H_{h,h'}}\right)\circ U_{g,g',h,h'}.\]
\end{itemize}
\end{defn}
This definition is equivalent to the one provided by \cite{fusion}.

Given a quiver $Q$, The path algebra $kQ$ consists of formal linear combinations of paths (including length $0$) in $Q$, with product given by the linear extension of concatenation (with $0$ resulting when paths are not compatible). We recall that an ideal $I\subseteq kQ$ is called \textit{admissible} if it does not contain any path of length $1$ or less, but contains all paths of length $\ge n$ for some $n$. The pair $(Q,I)$ is called a \textit{bound quiver}. There is a well known theorem of Gabriel \cite{quiver}, that relates basic algebras and bound quivers. Given this theorem, we will define a related bicatory of bounded quivers. The following definitions are necessary for this construction.

\begin{defn}
    Given quivers $Q_A$ and $Q_B$ with a quiver connection $(\Gamma, U)$, we introduce a new notation for paths. A path of type ${}_a(m,n)_b, m,n\geq 0$ consists of a path of length $m$ in $Q_A$ starting at vertex $a$ followed by an edge in $\Gamma$, followed by a path of length $n$ in $Q_B$ ending at vertex $b$. A path of type ${}_a(m,-)_b$ is a path of length $m$ in $Q_A$ from $a$ to $b$, and likewise a path of type ${}_a(-,n)_b$ is a path of length $n$ in $Q_B$ from $a$ to $b$. Path types are written adjacent for composition as follows: a path of type ${}_a(m,-)_b$ composed with a path of type ${}_b(n,-)_c$ would be expressed as ${}_a(m,-)_b{}_b(n,-)_c\sim {}_a(m+n,-)_c$. 
\end{defn}

Note that given a connection $(U,\Gamma)$ between quivers, we can ``iterate" $U$ to define isomorphisms:

\[U^{n}_{g,h}:\bigoplus E^{G}_{g,g_{1}}\otimes E^{G}_{g_1,g_{2}}\dots \otimes E^{G}_{g_{n-1},g_{n}} \otimes \Gamma_{g_{n},h}\to \bigoplus \Gamma_{g,h_{1}}\otimes E^H_{h_{1},h_{2}}\dots \otimes E^H_{h_{n},h},\] 

where on the left, the direct sum is over paths of type ${}_g(n,0)_h$, while on the right it is over paths of type ${}_g(0,n)_h$. This is defined as follows:

\[U_{g,h}^n=(U_{g,h_{n-1}}\otimes\id\otimes\dots\otimes \id)\circ(\id\otimes U_{g_1,h_{n-2}}\otimes\dots \otimes \id)\circ\dots\circ(\id\otimes \dots \id \otimes U_{g_{n-1},h})\]

In order to ``move'' edges across a connection between bound quivers, we need a notion of compatibility between the quiver ideals, which this definition provides.

\begin{defn}Given bound quivers $(Q_A,I_A)$ and $(Q_B,I_B)$, we say a quiver connection $(\Gamma, U)$ between them is \textit{ideally connected} if for each linear combination of paths of type ${}_a(m_i,-)_b$ in $I_A$ and each edge of type ${}_b(0,0)_c$ we have 
$$\sum_i c_iU^{m_i}({}_a(m_i,-)_b {}_b(0,0)_c)=\sum_i\sum_j d_j {}_a(0,0)_{b'} {}_{b'}(-,m_i)_c,$$
and likewise for $U^{-1}$, with $\sum_i\sum_j d_j {}_{b'}(-,m_i)_c\in I_B$.
\end{defn}

As mentioned above, every connected basic algebra is isomorphic to a quotient of its path algebra by an admissible ideal. With this is mind, we will describe a modification of the above 2-category:

\begin{defn}
We define the 2-category $\BdQuivCon$ whose
\begin{enumerate}
    \item 
    Objects are pairs $(Q, I)$ where $Q$ is a quiver and $I$ is an admissible ideal (these pairs are called \textit{bound quivers}).
    \item
    1-morphisms are quiver connections that are ideally connected
    \item
    2-morphisms are the same as in $\QuivCon$
    
\end{enumerate}
\end{defn}

Note that in the subcategory of bound quivers where $I=\rad^n$ (the path algebras obtained from these quivers are called \textit{truncated} path algebras), it is easy to see that every quiver connection is ideally connected.

\section{A 2-category of basic algebras}
Following the definition of $\BdQuivCon$, we then construct a 2-category of basic algebras with the goal of extending the quiver theorem to an equivalence of 2-categories. Recall that a $k$-algebra $A$ with a complete set $\{e_1,\hdots,e_n\}$ of primitive orthogonal idempotents is \textbf{basic} if $e_iA\not\cong e_jA$ for all $i\neq j$. We refer the reader to \cite{basic} for a comprehensive introduction to basic algebras, but recall some standard facts and introduce notation. Recall that the Jacobson radical $\rad(A)$ is nilpotent, i.e. the chain of ideals 
$$0=\rad^{m}(A)\subseteq \rad^{m-1}(A)\subseteq \rad(A)\subseteq A$$
stabilizes. We denote the quotient map 
$$\pi^{n}: \rad^{n-1}(A)\rightarrow \rad^{n-1}(A)/\rad^{n}(A).$$
Furthermore, since $\rad(A)$ is an ideal, for any subalgebras $B\subseteq A$, $\rad^{n}(A)$ is a $B$-$B$ sub-bimodule of $A$ as a $B$-$B$ bimodule.

Starting from the basic algebra $A$, in order to pin down the pair $(Q,I)$ up to isomorphism we need more data. First, we need to make a choice of a complete system of primitive orthogonal idempotents. This is equivalent to a splitting of the canonical short exact sequence

$$\begin{tikzcd}\label{ahh}
0\arrow{r} & \rad{A}\arrow{r} & A\arrow[swap]{r}{\pi^1} & A/\text{rad}(A)\arrow{r}\arrow[swap, bend right]{l}{\delta^{1}} & 0\\
\end{tikzcd}$$

Here and throughout this paper, we will denote the canonical quotient map\\ $\rad^{n-1}(A)\rightarrow \rad^{n-1}(A)/\rad^{n}(A)$ by $\pi^{n}$.

$$\begin{tikzcd}\label{ahh}
 A\arrow[swap]{r}{\pi^1} & A/\rad{A}\arrow[swap, bend right]{l}{\delta^{1}} \\
\end{tikzcd}$$

Since $A$ is basic, $A/\rad(A)\cong \mathbbm{k}^{n}$, so we have a complete system of primitive orthogonal idempotents given by $\{f_{i}:=\delta^{1}(e_{i})\}$, where $e_{i}$ are the minimal idempotents in $\mathbbm{k}^{n}$. Conversely, a choice of idempotents gives a splitting using the same formula. This data is sufficient to pin down the quiver $Q_{A}$ of the algebra $A$ (up to a natural notion of equivalence), but not the required admissible ideal. For this, note that the inclusion 
$$\delta^{1}: A/\rad(A)\hookrightarrow A$$
naturally makes the ideals $\rad(A)$ and $\rad^{2}(A)$ into $A/\rad(A)$ bimodules.  The extra data we need to specify a particular admissible ideal is a splitting of the short exact sequence of $A/\rad(A)$ bimodules

$$\begin{tikzcd}\label{ahh}
0\arrow{r} & \rad^{2}(A)\arrow{r} & \rad(A)\arrow[swap]{r}{\pi^2} & \rad(A)/\rad^{2}(A)\arrow{r}\arrow[swap, bend right]{l}{\delta^{2}} & 0\\
\end{tikzcd}$$

which we will abbreviate

$$\begin{tikzcd}\label{ahh}
 \rad(A)\arrow[swap]{r}{\pi^2} & \rad(A)/\rad^{2}(A)\arrow[swap, bend right]{l}{\delta^{2}} \\
\end{tikzcd}$$

Note that in both of the above situations, the $\pi$ maps are canonical quotients, the $\delta$ maps are additional data, and the condition on $\delta^{2}$ depends on the choice of $\delta^{1}$.

\begin{defn}\label{quivdat} If $A$ is a basic algebra, we define quiver data to be a choice of linear maps $\delta^{1}:A/\rad(A)\rightarrow A$ and $\delta^{2}: \rad(A)/\rad^{2}(A)\rightarrow \rad(A)$ satisfying

\begin{enumerate}
\item 
$\pi^{1}\circ \delta^{1}=id_{A/\rad{A}}$ and $\pi^{2}\circ \delta^{2}=id_{\rad(A)/\rad^{2}(A)}$
\item 
For $x,y\in A/\rad(A)$, $\delta^{1}(xy)= \delta^{1}(x)\delta^{1}(y)$
\item 
For $x,y\in A/\rad(A)$ and $z\in \rad(A)/\rad^{2}(A)$, $\delta^{2}(\delta^{1}(x)z\delta^{1}(y))=\delta^{1}(x)\delta^{2}(z)\delta^{1}(y)$

\end{enumerate}

\end{defn}

In the above, we are implicitly using the quotient $\delta^{1}(A/\rad(A))$ bimodule  on $\rad(A)/\rad^{2}(A)$.

Given quiver datum on a basic algebra, we have the following construction of a quiver:

\begin{defn}
Let $A$ be a basic finite dimensional $k$-algebra with choices $\delta^{1}_{A}$ and $\delta^{2}_{A}$ as above. The bound quiver of $A$, denoted $(Q_A, I_{A})$, is constructed as follows:
    \begin{enumerate}
        \item The vertex set $V(Q_A)$ is the set of minimal idempotents of the commutative semisimple alegbra $A/\rad(A)$.
        \item  Given vertices $e_{a},e_{b}\in V(Q_A)$, choose a basis $\{t_i\}$ for the $\mathbbm{k}$-vector space $$\delta^{1}_{A}(e_{a})(\delta^{2}_{A}(\rad(A)/\rad^{2}(A))\delta^{1}_{A}(e_b)$$ and define the edge set $E(e_{a}\rightarrow e_{b}):=\{t_i\}$.
        \item 
        The assignment of the edge $t_{i}$ to the corresponding element of $A$ defines a surjective homomorphism $\rho:kQ_{A}\rightarrow A$. Define $I_{A}:=\text{ker}(\rho)$.
        The pair $(Q_{A}, I_{A})$ is a bound quiver.
        
    \end{enumerate}

\end{defn}

Then, even given $\delta^{1}$ and $\delta^{2}$, the pair $(Q_{A}, I_{A})$ ostensiby depends on the choice of basis $\{t_{i}\}$ for $\delta^{1}_{A}(e_{a})(\delta^{2}_{A}(\rad(A)/\rad^{2}(A))\delta^{1}_{A}(e_b)$. But given any other basis $\{s_{i}\}$, the change of basis matrix gives a linear map $S_{a,b}$ between edge sets, which induce a graded isomorphism of the path algebras taking admissible ideal to admissible ideal.

\begin{rem} Usually in the statement and presentation of the above correspondence, a specific choice of $\delta^{2}_{A}$, is not given, and thus the ideal $I_{A}$ is \textit{not} canonically defined. We are  usually just satisfied with the existence portion of the theorem (choose some $\delta^{2}_{A}$). However, our goal is to extend the assignment $A\mapsto (Q_{A}, I_{A})$ into a 2-functor for certain 2-categories, which requires us to take the data $\delta^{1}_{A}$ and $\delta^{2}_{A}$ into account.
\end{rem}

Just as the semisimple quotient $A/\rad A$ is used to define the vertices of $Q_A$, we want a similar quotient to build the edges of the quiver connection corresponding to some $A-B$ bimodule $M$. The following definition enables us to ``semisimplify'' a bimodule:

\begin{defn}
Let $A,B$ be associative algebras and $M$ an $A$-$B$ bimodule. We say $M$ is radically symmetric if $\rad{A}M=M\rad{B}$, in which case we just refer to $\rad M$.
\end{defn}

If $M$ is radically symmetric, then $M/(\rad{A}M)=M/(M\rad{B})$ is an $A/\rad{A}-B/\rad{B}$ bimodule. We also note that $\rad^{2}(A)M=M\rad^{2}(B)$, so $\rad^{2}(M)$ is an unambiguously defined sub-bimodule. In particular we see $$\rad^{n}(A)\rad^{m}(M)=\rad^{n+m}(M)=\rad^{m}(M)\rad^{n}(B).$$ Clearly if $f:M\rightarrow N$ is morphism of radically symmetric $A$-$B$ bimodules, $f(\rad^{n}(M))\subseteq \rad^{n}(N)$.

We have a natural extension of the module multiplication map 

$$\rad^{n-1}(A)/\rad^{n}(A)\otimes \rad^{m-1}(M)/\rad^{m}(M)\rightarrow \rad^{n+m-2}(M)/\rad^{n+m-1}(M)$$ 

given by

$$(r+\rad^{n}(A))\otimes (s+\rad^{m}(M))\mapsto rs+\rad^{n+m-1}(M)$$

\noindent This is easily seen to be well-defined. A similar statement is true for the right $B$ action.

\begin{defn}\label{bimquivdat} If $(A,\delta^{1}_{A},\delta^{2}_{A})$ and $(B, \delta^{1}_{B}, \delta^{2}_{B})$ are basic algebras with quiver data, quiver data for a radically symmetric $A$-$B$ bimodule $M$ consists of $\delta^{1}_{M}: M/\rad(M)\rightarrow M$ and $\delta^{2}_{M}:\rad(M)/\rad^{2}(M)\rightarrow \rad(M)$ such that

\begin{enumerate}
\item 
$\pi^{1}_{M}\circ \delta^{1}_{M}=id_{M/\rad{M}}$ and $\pi^{2}_{M}\circ \delta^{2}_{M}=id_{\rad(M)/\rad^{2}(M)}$
\item 
For $x\in A/\rad(A)$, $y\in B/\rad(B)$, $u\in M/\rad(M)$ and $w\in \rad(M)/\rad^{2}(M)$,

$$\delta^{1}_{M}(\delta^{1}_{A}(x)u\delta^{1}_{B}(y))= \delta^{1}_{A}(x)\delta^{1}_{M}(u)\delta^{1}_{B}(y)$$

and 

$$\delta^{2}_{M}(\delta^{1}_{A}(x)w\delta^{1}_{B}(y))= \delta^{1}_{A}(x)\delta^{2}_{M}(w)\delta^{1}_{B}(y)$$

\item (compatibility)
For $x\in \rad(A)/\rad^{2}(A)$, $x\in \rad(B)/\rad^{2}(B)$, and $m\in M/\rad(M)$

$$\delta^{2}_{A}(x)\delta^{1}(m)=\delta^{2}(xm)$$

and

$$\delta^{1}(m)\delta^{2}_{B}(y)=\delta^{2}(my)$$

\end{enumerate}

\end{defn}

Note that we interpret $xm, my\in \rad{M}/\rad^{2}(M)$ as per the discussion before the definition.

\begin{defn}\label{bsa}
We define a 2-category $\BSA^{\circ}$ whose 

\begin{enumerate}
    \item 
    Objects are finite dimensional basic algebras with quiver data $(A, \delta^{1}_{A}, \delta^{2}_{A})$.
    \item
    1-morphisms are bimodules with quiver data $(M, \delta^{1}_{M}, \delta^{2}_{M})$
    \item
    2-morphisms $f:(M, \delta^{1}_{M}, \delta^{2}_{M})\rightarrow (N, \delta^{1}_{N}, \delta^{2}_{N})$ are $(M,N)$ bimodule intertwiners satisfying

    $$f\circ \delta^{1}_{M}=\delta^{1}_{N}\circ \widetilde{f}$$

    where $\widetilde{f}:M/\rad{M}\rightarrow N/\rad{N}$ is the induced intertwiner.

\end{enumerate}
\end{defn}

We are interested in the 2-category $\BSA$, which is a restriction of $\BSA^{\circ}$ to only dualizable 1-morphisms. The following lemma shows that given such a dualizable 1-morphism, we can recover a projective basis for the corresponding bimodule, which we later use to build the quiver connection.

\begin{lem}\label{lifted basis}
    Let $(M, \delta^{1}_{M}, \delta^{2}_{M})$ be an $A-B$ 1-morphism in $\BSA$. Then any standard basis for the lifted subspace $\delta^1_M(M/\rad M)$ is both a left projective basis and a right projective basis for $M$.
\end{lem}
\begin{proof}
If $(M, \delta^{1}_{M}, \delta^{2}_{M})$ is (right) dualizable then there exists a unit 2-morphism $$\eta:M^*\otimes_B M$$ satisfying the zigzag identities. Then $\eta$ is completely defined by the element $$\eta(1)=\sum\eta_i^*\otimes\eta_i\in M^*\otimes_A M.$$ Since the identity is contained in the subspace lifted by $\delta^1_B$, it follows that the $\eta_i$ are in the subspace lifted by $\delta^1_M$. The $\eta_i$ form a (right) projective basis for $M$, so it follows that we can choose a (right) projective basis for $M$ that lives in the desired subspace.

Now let $\{b_j\}$ be a standard basis for $\delta^1_M(M/\rad M)$, and $\{l_i\}$ be a left projective basis with dual projective basis $g_i:M_B\to B$ such that $\sum_i r_ig_i(m)=m$ for all $m\in M$. Set $l_i=\sum_j b_jL_{ji}$, $L_{ji}\in F$, and define $$f_j(\cdot):=\sum_i L_{ji}g_i(\cdot):M_B\to B^n.$$ Then 
\begin{align*}
    \sum_j b_jf_j(m)&=\sum_j b_j\sum_i L_{ji}g_i(m)\\
    &=\sum_{i,j}b_jL_{ji}g_i(m)\\
    &=\sum_il_ig_i(m)=m
\end{align*} for all $m\in M$, so it follows that $\{b_j\}$ is a left projective basis for $M$.

Similarly, let $\{r_i\}$ be a right projective basis in the lifted subspace with $h_i: \prescript{}{A}{M}\to A$ such that $\sum_i h_i(m)r_i=m$ for all $m\in M$. Then set $r_i=\sum_j R_{ji}(b_j), R_{ij}\in F$, and define $k_j(\cdot):=\sum_i R_{ji}h_i:\prescript{}{A}{M}\to A^n$. Then
\begin{align*}
    \sum_j k_j(m)b_j&=\sum_j \sum_i(R_{ji}h_i(m))b_j\\
    &=\sum_{i,j}R_{ji}h_i(m)b_j\\
    &=\sum_ih_i(m)r_i=m
\end{align*} for all $m\in M$, so it follows that $\{b_j\}$ is a right projective basis for $M$. Thus $\{b_j\}$ is a 2-sided projective basis for $M$ in the lifted subspace.

\end{proof}
This result shows that given a basis $\{b_j\}$ for $\delta^1_M(M/\rad M)$, any element $m\in M$ can be uniquely split into products $$m=\sum_jb_jf_j(m)=\sum k_j(n)b_j,$$ where $f_j(m)\in B$ and $k_j(m)\in A$.

The following lemmas are necessary for defining composition of 1-morphisms in $\BSA$ and are almost certainly well known to experts. 

\begin{lem}\label{Splittingradicalsovertensor}
    Let $M$ and $N$ be radically symmetric $A$-$B$ and $B$-$C$ bimodules respectively. Then there exists a canonical isomorphism of $A$-$C$ bimodules $f_{M,N}:M\otimes_B N/\rad M\otimes_B N\cong M/\rad M\otimes_B N/\rad N$.
\end{lem}
\begin{proof}
    Define $f_{M,N}:(M\otimes_B N)/(\rad M\otimes_B N)\to (M/\rad M)\otimes_B (N/\rad N)$ by considering the map 
    
    $$\tilde{f_{M,N}}: M\times N\rightarrow (M/\rad M)\otimes_B (N/\rad N),$$ 
    
    $$\tilde{f_{M,N}}(m,n):=(m+\rad N)\otimes_{B} (n+\rad M).$$
    
\noindent This is clearly $B$-balanced, hence uniquely extends to a bimodule map: $$f_{M,N}: M\otimes_{B} N\rightarrow M/\rad{M}\otimes_{B} N/\rad N.$$ We claim that $\rad (M\otimes_{B} N)\subseteq \text{ker}(f_{M,N})$. Note that $$\rad (M\otimes_{B} N)=(M\otimes_{B} N)\rad C=M\otimes_{B} \rad N,$$ and thus is spanned by simple tensors of the form $m\otimes_{B} n$ with $n\in\rad N$. Then 

\begin{align*}
    f_{M,N}(m\otimes n)&=(m+\rad M)\otimes (n+\rad N)\\
    &=(m+\rad M)\otimes 0+\rad N\\
    &=0.
\end{align*}

Thus $f_{M,N}$ descends to a well-defined bimodule map: $$f_{M,N}:(M\otimes_{B} N)/\rad (M\otimes_B N)\rightarrow M/\rad{M}\otimes_{B} N/\rad N.$$

Now define

$$f^{-1}_{M,N}:M/\rad{M}\otimes_{B} N/\rad N\rightarrow (M\otimes_{B} N)/\rad (M\otimes_B N)$$

by considering the universal extension of the $B$-balanced map 

$$\tilde{f}^{-1}_{M,N}: M/\rad{M}\times N/\rad N\rightarrow (M\otimes_{B} N)/\rad (M\otimes_B N)$$

$$\tilde{f}^{-1}_{M,N}(m+\rad M, n+\rad N):=m\otimes_{B}n +\rad (M\otimes_{B} N)$$

Note $\widetilde{f}^{-1}_{M,N}$ is well-defined since if $m-m^{\prime}=r\in \rad M$ and $n-n^{\prime}=s\in \rad N$, then

\begin{align*}
m\otimes_{B}n &= (m^{\prime}+r)\otimes_{B} (n^{\prime}+s)\\
&=m^{\prime}\otimes_{B} n^{\prime}+r\otimes_{B} n^{\prime}+m^{\prime}\otimes_{B} s+r\otimes_{B} s
\end{align*}

But by radical compatibility, $\rad M\otimes_{B} N=M\otimes_{B} \rad N=\rad (M\otimes_{B} N)$, so $r\otimes_{B} n^{\prime}+m^{\prime}\otimes_{B} s+r\otimes_{B} s\in M\otimes_{B} \rad N=\rad (M\otimes_{B} N)$. Therefore $\tilde{f}^{-1}_{M,N}$ extends uniquely to a linear (bimodule) map $f^{-1}_{M,N} $ as desired, which is clearly the inverse to $f_{M,N}$.

\end{proof}

\begin{lem}\label{Splittingradicalsquaredsovertensor}
    Let $M$ and $N$ be radically symmetric $A$-$B$ and $B$-$C$ bimodules respectively. Then there exists a canonical isomorphism of $A$-$C$ bimodules $f_{M,N}:\rad(M\otimes_B N)/\rad^2 (M\otimes_B N)\cong M/\rad M\otimes_B \rad N/\rad^2 N$.
\end{lem}

\begin{proof}
    Define $g_{M,N}:\rad(M\otimes_B N)/(\rad^2 M\otimes_B N)\to (M/\rad M)\otimes_B (\rad N/\rad^2 N)$ by considering the map 
    
    $$\tilde{g_{M,N}}: M\times \rad N\rightarrow (M/\rad M)\otimes_B (\rad N/\rad^2 N),$$ 
    
    $$\tilde{g_{M,N}}(m,n):=(m+\rad M)\otimes_{B} (n+\rad^2 N).$$

    \noindent This is clearly $B$-balanced, hence uniquely extends to a bimodule map $g_{M,N}: \rad(M\otimes_{B} N)\rightarrow M/\rad{M}\otimes_{B} \rad N/\rad^2 N$. We claim that $\rad^2 (M\otimes_{B} N)\subseteq \text{ker}(f_{M,N})$. Note that $\rad^2 (M\otimes_{B} N)=(M\otimes_{B} N)\rad^2 C=M\otimes_{B} \rad^2 N$, and thus is spanned by simple tensors of the form $m\otimes_{B} n$ with $n\in\rad^2 N$. Then 

$$g_{M,N}(m\otimes n)=(m+\rad M)\otimes (n+\rad^2 N)=(m+\rad M)\otimes 0+\rad^2 N=0$$

Thus $g_{M,N}$ descends to a well-defined bimodule map $$g_{M,N}:\rad(M\otimes_{B} N)/\rad^2 (M\otimes_B N)\rightarrow M/\rad{M}\otimes_{B} \rad N/\rad^2 N.$$

Now define

$$g^{-1}_{M,N}:M/\rad{M}\otimes_{B} \rad N/\rad^2 N\rightarrow \rad(M\otimes_{B} N)/\rad^2 (M\otimes_B N)$$

by considering the universal extension of the $B$-balanced map 

$$\tilde{g}^{-1}_{M,N}: M/\rad{M}\times \rad N/\rad^2 N\rightarrow \rad(M\otimes_{B} N)/\rad^2 (M\otimes_B N)$$

$$\tilde{g}^{-1}_{M,N}(m+\rad M, n+\rad^2 N):=m\otimes_{B}n +\rad^2 (M\otimes_{B} N)$$

Note $\widetilde{g}^{-1}_{M,N}$ is well-defined since if $m-m^{\prime}=r\in \rad M$ and $n-n^{\prime}=s\in \rad^2 N$, then

\begin{align*}
m\otimes_{B}n &= (m^{\prime}+r)\otimes_{B} (n^{\prime}+s)\\
&=m^{\prime}\otimes_{B} n^{\prime}+r\otimes_{B} n^{\prime}+m^{\prime}\otimes_{B} s+r\otimes_{B} s
\end{align*}

But by radical compatibility, $$\rad^2 M\otimes_{B} N=M\otimes_{B} \rad^2 N=\rad^2 (M\otimes_{B} N),$$ so $$r\otimes_{B} n^{\prime}+m^{\prime}\otimes_{B} s+r\otimes_{B} s\in M\otimes_{B} \rad^2 N=\rad^2 (M\otimes_{B} N).$$ Therefore $\tilde{g}^{-1}_{M,N}$ extends uniquely to a linear (bimodule) map $g^{-1}_{M,N} $ as desired, which is clearly the inverse to $g_{M,N}$.
\end{proof}

Note there is a natural surjective $A$-$C$ bimodule map $$\varphi:M\otimes_{B/\rad B} N\rightarrow M\otimes_{B} N.$$ Here, we are viewing $M$ and $N$ as right (left) $B/\rad B$ modules via the lifting\\
$B/\rad B\rightarrow B$ given by the choice of primitive orthogonal idempotents.

Now we can define the compositions of $1$ and $2$ morphisms in $\textbf{\BSA}$.

\begin{itemize}

\item
Composition of $1$-morphisms $(M,\delta^1_M,\delta^2_M)\in \textbf{\BSA}(A,B)$, and $(N,\delta^1_N,\delta^2_N)\in \textbf{\BSA}(B,C)$:
 
 $$(M,\delta^1_M,\delta^2_M)\otimes(N,\delta^1_N,\delta^2_N):=(M\otimes_{B} N,\delta^1_M\boxtimes\delta^1_N,\delta^1_M\boxtimes\delta^2_N),$$ where $\delta^1_M\boxtimes\delta^1_N:(M\otimes_B N)/\rad (M\otimes_B N)\to M\otimes_B N$ is defined as
 $$\varphi\circ(\delta^1_M\otimes_{B/\rad B}\delta^1_N)\circ f_{M,N},$$
 where $f_{M,N}$ is the isomorphism from Lemma \ref{Splittingradicalsovertensor}, and
$\delta^2_M\boxtimes\delta^2_N:(M\otimes_B N)/\rad (M\otimes_B N)\to M\otimes_B N$ is defined as
$$\varphi\circ(\delta^1_M\otimes_{B/\rad B}\delta^2_N)\circ g_{M,N},$$
where $g_{M,N}$ is the  isomorphism from Lemma $\ref{Splittingradicalsquaredsovertensor}$.

\item Vertical composition of 2-morphisms $g:(M,\delta^1_M,\delta^2_M)\to(N,\delta^{1}_{N},\delta^{2}_{N} )\in\BSA(A,B)$ and $h:(N,\delta^1_N,\delta^2_N)\to (P,\delta^1_P,\delta^2_P)$ is vertical composition of bimodule intertwiners.
\item Horizontal composition of 2-morphisms $g:(M,\delta^1_M,\delta^2_M)\to (N,\delta^1_N,\delta^2_N)$ and $h:(X,\delta^1_X,\delta^2_X)\to (Y,\delta^1_Y,\delta^2_Y)$ similarly is the horizontal composition of bimodule intertwiners.

 \end{itemize}

 \section{An equivalence of 2-categories}
 In this section we establish an equivalence between 2-categories $\BdQuivCon$ and $\BSA$. 

 First construct maps for a functor $\mathcal{P}': \QuivCon \to \PathAlg$. Given a quiver $Q$ with vertices $v_i$ and edges $e_j$, define $\mathcal{P'}(Q)=(kQ,\delta^1_{kQ},\delta^2_{kQ})$, where $\delta^1_{kQ}$ is defined on the basis $\{v_i+\rad kQ\}$ by $(v_i+\rad kQ)=v_i$ and $\delta^2_{kQ}$ is defined on the basis $\{e_j+\rad^2 kQ\}$ by $(e_j+\rad^2 kQ)=e_j$. It is clear that these maps satisfy the conditions for quiver data.

Given a $Q_A-Q_B$ quiver connection $(U,\Gamma)$, define $\mathcal{P}'(U,\Gamma)$ to be the right $kQ_B$ module $M=k(\Gamma Q_B)$ with paths of type ${}_a(0,n)_b$ forming a basis and the right action being concatenation of paths.  Then the set of edges $\{\gamma_i\}\in\Gamma$ form a right projective basis, and functionals $r_i$, defined by $$r_i(\sum_j \gamma_j b_j)=b_i,$$ form a right dual projective basis.
Define the left action of $kQ_A$ on $M$ by first using the isomorphisms $U^n$ to send paths to $kQ_B$, then act on the right by concatenation. We can see $\{\gamma_i\}$ is also a left projective basis, with functionals $l_i$ as a left dual projective basis, defined as follows: let $$\gamma b=\sum_j \gamma_j b_j$$ and $$\sum_j((U^n)^{-1}(\gamma_j b_j)=\sum_j a_j\gamma_j.$$ Then $l_i(\gamma b)=a_i$. Thus $P'(U,\Gamma)$ is a left and right projective module with the edges in $\Gamma$ forming a two-sided projective basis. Define 
$$\delta^1_{M}(\gamma_i+\rad M)=\gamma_i$$
and 
$$\delta^2_{M}(\gamma_ie_j+\rad^2 M)=\gamma_ie_j$$
for any edge $e_j\in Q_B$. These maps clearly satisfy the conditions for bimodule quiver data. This gives a 1-morphism in $\PathAlg$.

Lastly, given a 2-morphism $f=\{f_{i,j}\}$ between quiver connections $(U,\Gamma)$ and $(V,\Delta)$ define $\mathcal{P}'(f)$ to be the bimodule map $f^{\circ}$ defined by $$f^{\circ}(\sum \gamma b)=\sum f_{i,j}(\gamma) b.$$ This is clearly a bimodule intertwiner between $k\Gamma Q_B$ and $k\Delta Q_B$. 

Now we can consider all of these maps under quotients by bound ideals and define the functor $\mathcal{P}:\BdQuivCon\to\BSA$:
\begin{enumerate}
    \item Given a bound quiver $(Q, I)$, define $$\mathcal{P}(Q,I)=(kQ/I,\delta^1_{kQ/I},\delta^2_{kQ/I})$$
    \item Given an ideally compatible $(Q_A,I_A)-(Q_B,I_B)$ connection $(U,\Gamma)$, define $$\mathcal{P}(U,\Gamma)=(M/MI_B, \delta^1_{M/MI_B},\delta^2_{M/MI_B}),$$ which we will abuse notation to express as $(M/I_B, \delta^1_{M/I_B},\delta^2_{M/I_B})$
    \item Given a 2-morphism $f=\{f_{i,j}\}$ between quiver connections $(U,\Gamma)$ and $(V,\Delta)$ define $\mathcal{P}(f)$ to be the bimodule map $f^{\circ}$ defined by $$f^{\circ}(\sum \gamma b)=\sum f_{i,j}(\gamma) b$$
\end{enumerate}

We next construct a natural transformation $\mu_{\Gamma,\Delta}:\mathcal{P}(\Gamma\otimes\Delta)\to\mathcal{P}(\Gamma)\otimes\mathcal{P}(\Delta)$ using a similar argument to lemmas $\ref{Splittingradicalsovertensor}$ and $\ref{Splittingradicalsquaredsovertensor}$.

\begin{lem}\label{Splittingidealsovertensor}
    Let $(Q_A,I_A), (Q_B,I_B),$ and $ (Q_C,I_C)$ be bound quivers, and let $(U,\Gamma)$ and $(V,\Delta)$ be ideally connected $(Q_A,I_A)-(Q_B,I_B)$ and $(Q_B,I_B)-(Q_C,I_C)$ quiver connections respectively. Then there exists a natural, canonical isomorphism of $kQ_A/I_A$-$kQ_C/I_C$ bimodules $\mu_{\Gamma,\Delta}:k(\Gamma\otimes_{Q_B} \Delta)/ I_C \cong k\Gamma/ I_B\otimes_{kQ_B/I_B} k\Delta/I_C$.
\end{lem}
\begin{proof}
    Define $\mu_{\Gamma,\Delta}:k\Gamma\otimes_{Q_B} \Delta/I_C\to k\Gamma/I_B\otimes_{kQ_B/I_B} k\Delta/I_C$ by considering the map 
    
    $$\tilde{\mu_{\Gamma,\Delta}}: k(\Gamma\times\Delta)\rightarrow k\Gamma/I_B\otimes_{kQ_B/I_B} k\Delta/I_C,$$ 
    
    $$\tilde{\mu_{\Gamma,\Delta}}(\gamma,\delta):=(\gamma+I_B)\otimes_{kQ_B/I_B} (\delta+I_C).$$
    
\noindent This is clearly $Q_B$-balanced, hence uniquely extends to a bimodule map $$\mu_{\Gamma,\Delta}: k(\Gamma\otimes_{Q_B}\Delta)\rightarrow k\Gamma/I_B\otimes_{kQ_B/I_B} k\Delta/I_C.$$ We claim that $k(\Gamma\otimes\Delta)I_C\subseteq \text{ker}(\mu_{\Gamma,\Delta})$. Note that since the quiver connections are ideally connected, it follows that $k(\Gamma\otimes\Delta)I_C$ is spanned by simple tensors of the form $\delta\otimes \gamma$ with $\gamma\in k(\Gamma\otimes\Delta)I_C$. Then 

$$\mu_{\Gamma,\Delta}(\gamma\otimes \delta)=(\gamma+I_B)\otimes (\delta+I_C)=(m+\rad M)\otimes 0+I_C=0$$

Thus $\mu_{\Gamma,\Delta}$ descends to a well-defined bimodule map $\mu_{\Gamma,\Delta}:k\Gamma\otimes_{Q_B} \Delta/I_C \cong k\Gamma/I_B\otimes_{kQ_B/I_B} k\Delta/I_C$.

Now define

$$\mu^{-1}_{\Gamma,\Delta}:k\Gamma/I_B\otimes_{kQ_B/I_B} k\Delta/I_C\rightarrow k\Gamma\otimes_{Q_B} \Delta/I_C$$

by considering the universal extension of the $kQ_B/I_B$-balanced map 

$$\tilde{\mu}^{-1}_{\Gamma,\Delta}: k\Gamma/I_B\times k\Delta/I_C\rightarrow k(\Gamma\otimes_{Q_B} \Delta)/I_C$$

$$\tilde{\mu}^{-1}_{\Gamma,\Delta}(\gamma+I_B, \delta+I_C):=\gamma\otimes_{B}\delta +k(\Gamma\otimes\Delta)I_C$$

Note $\widetilde{\mu}^{-1}_{\Gamma,\Delta}$ is well-defined since if $\gamma-\gamma^{\prime}=r\in I_B$ and $\delta-\delta^{\prime}= I_C$, then

\begin{align*}
\gamma\otimes_{kQ_B/I_B}\delta &= (\gamma^{\prime}+r)\otimes_{Q_B} (\delta^{\prime}+s)\\
&=\gamma^{\prime}\otimes_{Q_B} \delta^{\prime}+r\otimes_{Q_B} \delta^{\prime}+\gamma^{\prime}\otimes_{Q_B} s+r\otimes_{Q_B} s
\end{align*}

But by ideal connectivity $r\otimes_{Q_B} \delta^{\prime}+\gamma^{\prime}\otimes_{Q_B} s+r\otimes_{Q_B} s\in k(\Gamma\otimes_{Q_B} \Delta)I_C$. Therefore $\tilde{\mu}^{-1}_{\Gamma,\Delta}$ extends uniquely to a linear (bimodule) map $\mu^{-1}_{\Gamma,\Delta} $ as desired, which is clearly the inverse to $\mu_{\Gamma,\Delta}$.

Now we show naturality. Let $(U',\Gamma')$ and $(V',\Delta')$ be ideally connected $(Q_A,I_A)-(Q_B,I_B)$ and $(Q_B,I_B)-(Q_C,I_C)$ quiver connections respectively, let $f$ be a 2-morphism between $(U,\Gamma)$ and $(U',\Gamma')$, and let $g$ be a 2-morphism between $(V,\Delta)$ and $(V',\Delta')$. We must show that $$\mathcal{P}(f)\otimes\mathcal{P}(g)\circ \mu_{\Gamma\otimes\Delta}=\mu_{\Gamma'\otimes\Delta'}\circ \mathcal{P}(f\otimes g).$$ Note that $P(\Gamma\otimes\Delta)$ is spanned by simple tensors of the form $\gamma\otimes\delta c$ where $\gamma,\delta$ are in the image of $\delta^1_{k\Gamma},\delta^1_{k\Delta}$ respectively and $c\in kQ_C/I_C$. Taking such a simple tensor, we have $$\mathcal{P}(f)\otimes\mathcal{P}(g)\circ \mu_{\Gamma\otimes\Delta}(\gamma\otimes\delta c)=\sum f_{g,h}(\gamma)\otimes\sum f_{h,k}(\delta)c$$ and $$\mathcal{P}(f\otimes g)(\gamma\otimes\delta c)=\sum_{i,k} f\otimes  g_{i,k}(\gamma\otimes\delta)c.$$ By definition this is equal to $\sum_{i,j,k}(f_{i,j}(\gamma)\otimes g_{j,k}(\delta))c$, so it follows that 
\begin{align*}
    \mu_{\Gamma'\otimes\Delta'}\circ \mathcal{P}(\gamma\otimes\delta)&=\mu_{\Gamma'\otimes\Delta'}(\sum_{i,j,k}(f_{i,j}(\gamma)\otimes g_{j,k}(\delta))c)\\
    &=\sum f_{g,h}(\gamma)\otimes\sum f_{h,k}(\delta)c,
\end{align*}as desired. Thus $\mu$ is natural.
\end{proof}

It remains to show that $\mathcal{P}$ is an equivalence. Essential surjectivity on objects is provided by the quiver theorem along with a choice of quiver data.

\begin{lem}
    $\mathcal{P}$ is essentially surjective.
\end{lem}
\begin{proof}
    Let $(A,\delta^1_A,\delta^2_A)$ be a basic algebra with quiver data, and $Q_A$ by the quiver associated to $A$. By Gabriel's quiver theorem there exists an algebra homomorphism $\varphi:kQ_A\to A$ that is surjective with kernel $I=\ker\varphi$ that is an admissible ideal of $kQ_A$. Furthermore, the choice of quiver data parameterizes $\varphi$. By construction, $\mathcal{P}$ generates quiver data ($\delta^1_{kQ_A/I}, \delta^2_{kQ_A/I}$) that is compatible with the quiver data on $A$. Thus $\mathcal{P}$ is essentially surjective on equivalence classes of objects.
\end{proof}

The following lemma shows that given a left projective basis, a bimodule $M$ can be decomposed into a tensor product of the space spanned by the projective basis and the algebra that acts from the right.

\begin{lem}
    Let $(M,\delta^1_M,\delta^2_M)\in\BSA$. Then we can construct a unique bimodule isomorphism $M\cong \delta^1_M(M/\rad M)\otimes_{\delta^1_B(B/\rad B)} B$.
\end{lem}
\begin{proof}
    Recall that the quiver datum $\delta^1_M$ corresponds to a standard basis $\{m_j\}$ for $\delta^1_M(M/\rad M)$. By lemma \ref{lifted basis}, this is a left projective basis with dual projective basis $n_j$, so for any $m\in M$ we can uniquely write the sum $m=\sum_j m_j n_j(m)$. Since $M/\rad M$ decomposes into the direct sum $\bigoplus (M/\rad M)f_j$, choose a basis $\{m_j\}$ such that for all $i$, $m_i\in\delta^1_M((M/\rad M)f_i)$. Then define $g_M:M\to \delta^1_M(M/\rad M)\otimes_{\delta^1_B(B/\rad B)} B$ as $g_M(m)=\sum_j m_j\otimes n_j(m)$, and define $g^{-1}_M:\delta^1_M(M/\rad M)\otimes_{\delta^1_B(B/\rad B)} B\to M$ by $g^{-1}_M(m\otimes b)=mb$. Note that since $m\in\delta^1_M(M/\rad M)=\sum_j m_j$, it follows that $mb=\sum_j (m_jf_j)b=0$ if and only if $\sum_j m_j\otimes_{\delta^1_B(B)}b=0$. Thus $g_M^{-1}$ is well-defined and injective. It is clear that $g_M^{-1}$ is surjective, and thus that $g_M$ is its inverse.
\end{proof}

\begin{lem}\label{isomorphism}
Let $(A,\delta^1_A,\delta^2_A),(B,\delta^1_B,\delta^2_B)\in \BSA$ and $(M,\delta^1_M,\delta^2_M)\in \BSA(A,B)$. Then using the quiver data we can uniquely construct a family of linear isomorphisms
\begin{align*}
    U_{i,k}=:&\bigoplus_j e_i(\rad A/\rad^2 A)e_j\otimes_{A/\rad A} e_j(M/\rad M)f_k\to\\
    &\bigoplus_h e_i(M/\rad M )f_h\otimes_{B/\rad B} f_h(\rad B/\rad^2 B)f_k
\end{align*}
compatible with the left and right actions on $M$; if $$U_{i,j}(a+\rad^2 A\otimes m+\rad M)=\sum m'+\rad M\otimes b+\rad^2 B$$ then $$\delta^2_A(a+\rad A)\delta^1_M(m+\rad M)=\sum\delta^1_M(m'+\rad M)\delta^2_B(b+\rad^2 B).$$
Furthermore, the corresponding quiver connection is ideally connected; if $$\sum_s \delta^2_A(a_{s_1})\delta^2_A(a_{s_2})\hdots \delta^2_A(a_{s_t})\delta^1_M(m_s)=0$$ then $$\sum_sU^{s_t}_{i_s,j_s}(a_{s_1}\otimes a_{s_2}\otimes\hdots\otimes a_{s_t}\otimes m_s)=0.$$ 
\end{lem}
\begin{proof}
Let $a+\rad^2 A\in e_i(\rad A/\rad^2 A)e_j$ and $m+\rad M\in e_j(M/\rad M)f_k$. Note that due to the radical symmetry of $M$ we have that $\rad AM=M\rad B$, and by extension $(\rad A)M/\rad^2M=M\rad B/\rad^2M$. Since any basis for $M/\rad M$ lifts to a projective basis for $M$, and since $\delta^1_M$ determines a basis for $M/\rad M$, there exists unique $m_h+\rad M\in e_i(M/\rad M)f_h$ and $b_h+\rad^2 B\in f_h(\rad B/\rad^2 B)f_k$ such that
\begin{align*}
    \delta^2_A(a+\rad^2 A)\delta^1_M(m+\rad M)&=\delta^2_M(am+\rad^2M)\\
    &=\sum_h\delta^2_M( m_hb_h+\rad^2 M)\\
    &=\sum_h\delta^1_M(m_h+\rad M)\delta^2_B(b_h+\rad^2 B),
\end{align*}
so we define $$U_{i,k}(a+\rad^2 A\otimes m+\rad M)=\sum_h (m_h+\rad M)\otimes (b_h+\rad^2 B).$$ Since any basis for $M/\rad M$ lifts to a 2-sided projective basis for $M$, we can similarly define $U^{-1}_{i,k}$, and it is clear that both of these maps are well-defined.

Similar to our construction in $\BdQuivCon$, these isomorphisms can be iterated:
\[U_{i,j}^n=:(U_{i,j_{n-1}}\otimes\id\otimes\dots\otimes \id)\circ(\id\otimes U_{i_1,j_{n-2}}\otimes\dots \otimes \id)\circ\dots\circ(\id\otimes \dots \id \otimes U_{i_{n-1},j})\]

By definition of the lifting maps, we have:
    \begin{align*}
        0&=\sum_s \delta^2_A(a_{s_1})\delta^2_A(a_{s_2})\hdots \delta^2_A(a_{s_t})\delta^1_M(m_s)\\
        &=\sum_s \delta^2_A(a_{s_1})\delta^2_A(a_{s_2})\hdots \delta^2_M(a_{s_t}m_s)\\
        &=\sum_s \delta^2_A(a_{s_1})\delta^2_A(a_{s_2})\hdots \delta^2_A(a_{s_t-1})\sum_{s'}\delta^1_M(m_{s'})\delta^2_B(b_{s'})\\
        &\hspace{3cm}\vdots\\
        &=\sum_{s'}\delta^1_M(m_{s'})\delta^2_B(b_{s_1})\delta^2_B(b_{s_2'})\hdots \delta^2_B(b_{s_t'}),
    \end{align*}
    and
    \begin{align*}
        &\sum_sU^{s_t}_{i_s,j_s}(a_{s_1}\otimes a_{s_2}\otimes\hdots\otimes a_{s_t}\otimes m_s)=\sum_{s'}m_s'\otimes b_{s_1'}\otimes b_{s_2'}\hdots\otimes b_{s_t'}=0,
    \end{align*}
    so the connection is ideally connected, as desired.
\end{proof}

Now we use the explicit isomorphisms between $\Gamma$ and $\delta^1_M(M/\rad M)$ and between $A$ and $kQ_A/I$ to construct a 2-morphism that is an isomorphism.

\begin{lem}
    $\mathcal{P}$ is essentially full on 1-cells.
\end{lem}
\begin{proof}
    Let $(A,\delta^1_A,\delta^2_A), (B,\delta^1_B,\delta^2_B)\in \BSA$  and $(M,\delta^1_M,\delta^2_M)\in \BSA(A,B)$. Construct a $(Q_A, I_A)-(Q_B,I_B)$ quiver connection ($\Gamma,U)$ where $\Gamma_{i,j}\cong e_i(M/\rad M)f_j$ as $k$-vector spaces and $U$ is constructed as in lemma \ref{isomorphism}. Then 
    $$\mathcal{P}(\Gamma, U)=(k\Gamma Q_B/I_B, \delta^1_{\Gamma Q_B/I_B}, \delta^2_{\Gamma Q_B/I_B}).$$
    It remains to show that there is a 2-morphism
    $$h:(M,\delta^1_M,\delta^2_M)\to (k\Gamma Q_B/I_B, \delta^1_{\Gamma Q_B/I_B}, \delta^2_{\Gamma Q_B/I_B})$$
    which is an isomorphism. By construction, there are isomorphisms $$\delta_{i,j}:(e_i(M/\rad M)f_j)\to \Gamma_{i,j}$$ and $$\varphi:B\to kQ_B/I_B.$$ Then using bimodule isomorphisms $g_M$ and $g^{-1}_{k\Gamma Q_B/I_B}$, we can build the bimodule isomorphism $$g^{-1}_{k\Gamma Q_B/I_B}\circ (\bigoplus(\delta_{ij})\otimes \Delta)\circ g_M:M\to k\Gamma Q_B/I_B,$$
    as desired.
\end{proof}

It is clear that $\mathcal{P}$ is fully faithful on 2-cells, and thus fully faithful. Since $\mathcal{P}$ is essentially surjective, $\mathcal{P}$ is an equivalence of 2-categories, completing the proof of Theorem (A).

\end{document}